\documentclass[11pt]{article}
\newtheorem{conjecture}{Conjecture}

\newtheorem{corollary}[conjecture]{Corollary}
\newtheorem{definition}[conjecture]{Definition}

\newtheorem{lemma}[conjecture]{Lemma}

\newtheorem{proposition}[conjecture]{Proposition}

\newtheorem{theorem}[conjecture]{Theorem}

\newenvironment{proof}{\noindent {\bf Proof} \hspace{.1cm}}{\hfill ${\bf QED}$ \\ \vspace{.15cm}}

\title{A simple criterion for non-relative hyperbolicity and one-endedness of groups}
\author{James W. Anderson, Javier Aramayona and Kenneth J. Shackleton\footnote{Supported by an EPSRC studentship}}

\date{23 November, 2005}

\begin{document}

\maketitle

\begin{abstract}
\noindent We give a combinatorial criterion that implies both  the
non-strong relative hyperbolicity and the one-endedness of a finitely generated group. We use
this to show that many important classes of groups do not admit a strong
relatively hyperbolic group structure and have one end. Applications include
surface mapping class groups, the Torelli group, (special)
automorphism and outer automorphism groups of most free groups, and the three-dimensional Heisenberg group. Our final application is to Thompson's group $F$.

\medskip

\noindent MSC 20F67 (primary), 20F65 (secondary)
\end{abstract}


\section{Introduction}

In recent years, the notion of relative hyperbolicity has become a
powerful method for establishing analytic and geometric properties
of groups, see for example Dadarlat and Guentner \cite{dadarlat-guentner},
Dahmani \cite{dahmani}, Osin \cite{osin-asymptotic}, \cite{osin-relatively},
Ozawa \cite{ozawa-boundary}, Yaman \cite{yaman}.  Relatively hyperbolic
groups, first introduced by Gromov \cite{gromov-hyperbolic} and
then elaborated on by various authors (see Farb \cite{farb-relatively},
Szczepa\'nski \cite{szczepanski}, Bowditch \cite{bowditch-relatively}), provide a natural
generalization of hyperbolic groups and geometrically finite
Kleinian groups.

When a finitely generated group $G$ is strongly hyperbolic relative to a
finite collection $L_1, L_2, \ldots, L_p$ of proper subgroups, it is often possible to deduce that $G$ has a given property provided the subgroups $L_j$ have the same
property. Examples of such properties include finite asymptotic dimension, see Osin \cite{osin-asymptotic}, exactness, see Ozawa \cite{ozawa-boundary}, and
uniform embeddability in Hilbert space, see Dadarlat and Guentner 
\cite{dadarlat-guentner}. In light of this, identifying a
strong relatively hyperbolic group structure for a given
group $G$, or indeed deciding whether or not one can exist,
becomes an important objective. 

One of the main
results of this note, Theorem \ref{our-main-theorem} in Section \ref{comm-graph}, asserts
that such a structure cannot exist whenever the group $G$
satisfies a simple combinatorial property, namely that its 
{\em commutativity graph} with respect to some generating set $S$ is
connected. We describe this graph in Section $3$. 

The other main result of this note, Theorem \ref{main-end-theorem} in Section \ref{ends-groups}, asserts that whenever $G$ has connected commutativity graph with respect to some set of generators, it has one end.

Using our main results,
we deduce that many well-known groups do not admit a
strong relatively hyperbolic group structure and have one end. These examples
include all but finitely many surface mapping class groups, the Torelli group of a
closed surface of genus at least 3, the (special) automorphism and outer automorphism groups of almost all free groups and the three-dimensional Heisenberg group. We remark that the one-endedness of the surface mapping class groups and  (special) automorphism and outer automorphism groups of free groups considered herein was previously established by Culler and Vogtmann  \cite{culler-vogtmann} using Bass-Serre theory. In Section $4.5$, we prove that Thompson's group $F$ has one end and is not strongly relatively hyperbolic using a minor variation of our main argument.

During the preparation of this work we were informed that
Behrstock, Dru\c{t}u, and Mosher \cite{bdm} have found an alternative argument
for the non-strong relative hyperbolicity of some of the groups treated here, using asymptotic cones and the description of
relative hyperbolicity due to Dru\c{t}u and Sapir \cite{drutu-sapir}.

\noindent 
{\bf{Acknowledgements.}} Parts of this paper were written during visits of the second author to the Bernoulli Center, Lausanne, and the third author to the Tokyo Institute of Technology. Both authors wish to express their gratitude to these institutions.  We would also like to thank Martin Dunwoody and Ian Leary for helpful conversations.

\section{Relatively hyperbolic groups}
\label{rel-hyp-group}

We assume that all groups appearing in this note are infinite, unless otherwise explicitly stated.

There are two related but inequivalent definitions of relative hyperbolicity that are commonly used, one due to Farb 
\cite{farb-relatively}, and the other developed by  Gromov \cite{gromov-hyperbolic}, Szczepa\'nski \cite{szczepanski}, and Bowditch \cite{bowditch-relatively}. As we only go into as much detail as is required for us to state our results, we refer 
the interested reader to the cited papers for a more extensive treatment.

We first give the definition given by Farb \cite{farb-relatively} and refer to this as {\em weak relative hyperbolicity}. For a group $G$, a finite generating set $S$ and a finite family of proper finitely generated subgroups $\{L_1, L_2, \ldots, L_m\}$, we 
form an augmentation ${\rm Cay}^{*}(G, S)$ of the Cayley graph ${\rm Cay}(G,S)$ as follows:  Give ${\rm Cay}(G, S)$ the 
path-metric obtained by declaring each edge to have length one.  Then, for each $1\le j\le m$, adjoin to ${\rm Cay}(G, S)$ a new vertex $v_{gL_j}$ for each coset $gL_j$ of $L_j$ and declare the distance between each new vertex $v_{g L_j}$ and each 
vertex in the associated coset $g L_j$ to be one.  We say that $G$ is {\em weakly hyperbolic relative to $L_1, L_2, 
\ldots, L_m$} if the resulting metric on ${\rm Cay}^{*}(G, S)$ is hyperbolic in the sense of Gromov. Farb \cite{farb-relatively} shows this definition does not depend on the choice of generating set. 

In the same paper, Farb introduces the notion of {\em bounded coset penetration} (BCP), which is a weak local finiteness 
property satisfied by many important examples of weakly relatively hyperbolic groups.  Roughly speaking, BCP imposes certain fellow-travelling conditions on 
pairs of quasi-geodesics on ${\rm Cay}(G, S)$ with the same endpoints that enter cosets of the subgroups $L_j$, $1\le j\le m$.

Bowditch \cite{bowditch-relatively} gives two equivalent dynamical notions of relative hyperbolicity, of which we recall 
the second. We will refer to this notion as {\em strong relative hyperbolicity}.  We say that a group $G$ is 
{\em strongly hyperbolic relative to the family $L_1, L_2, \ldots, L_m$ of proper finitely generated subgroups} if $G$ 
admits an action on a connected, hyperbolic graph $\mathcal{G}$ such that $\mathcal{G}$ is {\em fine} (that is, for each $n\in{\bf{N}}$, 
each edge of $\mathcal{G}$ belongs to only finitely many circuits of length $n$), there are only finitely many $G$-orbits of edges, each edge 
stabiliser is finite, and the stabilisers of vertices of infinite valence are precisely the conjugates of the 
$L_j$.

We note that strong relative hyperbolicity (with respect to some finite collection of proper finitely generated subgroups) is equivalent to weak relative hyperbolicity plus BCP (with respect to the same collection of subgroups), see Szczepa\'nski \cite{szczepanski} and Dahmani \cite{dahmani}. The BCP property is crucial: as noted in Szczepa\'nski \cite{szczepanski}, the group $\bf{Z} \oplus \bf{Z}$ is weakly, but not strongly, hyperbolic relative to the diagonal subgroup $\{ (m,m)\: |\: m\in \bf{Z}\}$. 


We mention here that a characterisation of strong relative hyperbolicity in terms of relative presentations and isoperimetric inequalities is given in Osin \cite{osin-relatively}.

\section{The commutativity graph}
\label{comm-graph}

We begin by defining a graph which attempts to capture the notion that a group is well generated by large abelian subgroups.

\begin{definition}[Commutativity graph]
Let $G$ be a group and let $S$ be a (possibly infinite) generating set for $G$, all of whose elements have infinite order. The {\em commutativity graph} $K(G, S)$ for $G$ with respect to $S$ is the simplicial graph whose vertex set is $S$ and in which distinct vertices $s$, $s'$ are connected by an edge if and only if there are non-zero integers $n_s$, $n_{s'}$ so that $\langle s^{n_s}, (s')^{n_{s'}}\rangle$ is abelian.  
\end{definition}

As long as there is no danger of confusion, we will use the same notation for elements of $S$ and vertices of $K(G,S)$.

Notice that for any $g\in G$, we have that $K(G,S)$ is connected if and only if $K(G, gSg^{-1})$ is connected.  Typically we shall only consider commutativity graphs in which adjacent vertices, rather than powers of adjacent vertices, commute.


Our main result about non-strong relative hyperbolicity may be stated as follows.  Recall that the {\em rank} of a finitely generated abelian group $A$ is the rank of some (and hence every) free abelian subgroup $A_0$ of finite index in $A$.


\begin{theorem} Let $G$ be a finitely generated group. Suppose there exists a (possibly infinite) generating set $S$ of cardinality at least two such that every element of $S$ has infinite order and $K(G, S)$ is connected. Suppose further that there exist adjacent vertices $s$, $s'$ of $K(G,S)$ and non-zero integers $n_s$, $n_{s'}$ so that $\langle s^{n_s}, (s')^{n_{s'}}\rangle$ is rank $2$ abelian.  
Then, $G$ is not strongly hyperbolic relative to any finite collection of proper finitely generated subgroups. 
\label{our-main-theorem}
\end{theorem}

The remainder of this section is dedicated to the proof of Theorem
\ref{our-main-theorem}. The main tool we use is the following
theorem on virtual malnormality for strongly relatively hyperbolic groups, which is contained in the work of Farb
\cite{farb-relatively} and Bowditch \cite{bowditch-relatively}, and
is explicitly stated in Osin \cite{osin-relatively} (Theorems $1.4$ and
$1.5$). 


\begin{theorem} Let $G$ be a finitely generated group that is strongly hyperbolic relative to the proper finitely generated subgroups $L_1,\ldots, L_p$.  Then,
\begin{enumerate}
\item For any $g_1$, $g_2\in G$, the intersection $g_1 L_j g_1^{-1}\cap g_2 L_k g_2^{-1}$ is finite for
$1\le j\ne k\le p$.
\item For $1\le j\le m$, the intersection $L_j \cap g L_j g^{-1}$ is finite for any $g\not\in L_j$.
\end{enumerate}
\label{relatively-malnormal}
\end{theorem}

Note that this immediately implies that if $g\in G$ has infinite order and if $g^k$ lies in a conjugate $hL_j h^{-1}$ of some $L_j$, then $g$ lies in that same conjugate $hL_j h^{-1}$ of $L_j$, since the intersection $hL_j h^{-1} \cap ghL_j h^{-1}g^{-1}$ then contains $\langle g^k\rangle$, which is infinite.

We also need the following lemma, which follows directly
from Theorems 4.16 and 4.19 of Osin \cite{osin-relatively}. We note that while this lemma is implicit in the literature, it has never (to the best
of our knowledge) been stated in this form, and so we include it as it may
be of independent interest. 

\begin{lemma}
Let $G$ be a finitely generated group that is strongly hyperbolic
relative to the proper finitely generated subgroups $L_1, L_{2}, \ldots,
L_p$.  If $A$ is an abelian subgroup of $G$ of rank at least two,
then $A$ is contained in a conjugate of one of the $L_j$.
\label{abelian-lemma}
\end{lemma}

We are now ready to prove Theorem \ref{our-main-theorem}.

\begin{proof} [of Theorem \ref{our-main-theorem}]
Suppose for contradiction that $G$ is strongly hyperbolic relative to the finite collection
$L_{1}, L_{2}, \ldots, L_{p}$ of proper finitely generated
subgroups. We first show that no conjugate of any 
$L_j$ can contain a non-zero power of an element of $S$. So, suppose there is some $g\in G$, some $s_0\in S$, and some $k\ne 0$ so that $s_{0}^{k} \in g L_j g^{-1}$ for some $1\le j\le p$.  (Note that, by the comment immediately following Theorem \ref{relatively-malnormal}, this implies that $s_0\in gL_j g^{-1}$ as well.)  

Let $s_{1}$ be any vertex of $K(G, S)$ adjacent to
$s_{0}$. As there are non-zero integers $n_0$, $n_1$ so that $\langle s_{0}^{n_0}, s_{1}^{n_1}\rangle$ is abelian, we see that $\langle s_{0}^{n_0 k}\rangle \subseteq g L_j g^{-1} \cap s_{1}^{n_1}g L_j g^{-1} s_{1}^{-n_1}$. However, since the
subgroup $\langle s_{0}^{n_0 k}\rangle $ of $G$ is infinite, Theorem
\ref{relatively-malnormal} implies that $s_1^{n_1} \in g L_j g^{-1}$.  By the comment immediately following Theorem \ref{relatively-malnormal}, we see that $s_1\in gL_jg^{-1}$.

Now, let $s$ be any element of $S$. By the connectivity of $K(G, S)$,
there is a sequence of elements $s_{0}, s_{1}, \ldots, s_{n}=s$
of $S$ such that $s_{k-1}$ and $s_k$ are adjacent in $K(G,S)$ for each $1\le k\le n$. The argument we have given above implies that if $s_{k-1} \in g L_j g^{-1}$,
then $s_{k} \in g L_j g^{-1}$. In particular, we have that $s\in g L_j g^{-1}$. Since $G$ is generated by $S$, it follows that $G$
and $g L_j g^{-1}$ are equal, contradicting the assumption that
the subgroup $L_j$ is proper.  We conclude that if $G$ is strongly hyperbolic relative to $L_{1},
L_{2}, \ldots, L_{p}$ then no conjugate of any $L_{j}$ can contain a non-zero power
of an element of $S$. 

Now, by assumption there exist adjacent vertices $t$ and $t'$ of $K(G,S)$ for which there exist non-zero integers $n_t$, $n_{t'}$ so that $A =\langle t^{n_t}, (t')^{n_{t'}}\rangle$ is rank $2$ abelian.  By Lemma \ref{abelian-lemma}, we see that $A$ must lie in some conjugate of some $L_{j}$. In particular, this conjugate of $L_{j}$ contains a non-zero power of an element of $S$, and we have a final contradiction.
\end{proof}

\section{Ends of groups}
\label{ends-groups}

In this section, we use the same argument as that used in the proof of Theorem \ref{our-main-theorem} to show that groups with connected commutativity graph have one end.  
We refer to Lyndon and Schupp \cite{lyndon-schupp} for standard facts about amalgamated free products and HNN extensions.

We do not give here a complete and formal definition of ends of groups; for this, the interested reader is referred to Epstein \cite{epstein}, 
Stallings \cite{stallings}, and Dicks and Dunwoody \cite{dicks-dunwoody}.  Let $G$ be a finitely generated group (which we have assumed to be infinite).  
Say that $G$ has {\em one end} if for some (and hence for every) finite generating set $S$, the Cayley graph ${\rm Cay}(G,S)$ has the property that the complement 
${\rm Cay}(G,S) \setminus B_n(e)$ is connected, where $B_n(e)$ is the ball of radius $n$ about the identity $e$.  Say that $G$ has {\em two ends} if $G$ is virtually ${\bf Z}$. 

Say that a finitely generated group $G$ has {\em infinitely many ends} if and only if $G$ admits either a {\em non-trivial amalgamated free product} 
decomposition $G =A\ast_C B$ or a {\em non-trivial HNN} decomposition $G =A\ast_C$, where (in either case) $C$ is finite. Here, by a trivial amalgamated free
product decomposition, we mean a decomposition of the form $G =A\ast_C B$ where $C$ is finite and $[A:C] =[B:C] =2$, and by a trivial HNN extension we mean an HNN extension of the form $G =C\ast_C$ where $C$ is finite. Note that, in 
these two cases, $C$ is a finite normal subgroup of $G$ and $G/C$ is infinite cyclic, and it follows $G$ has two ends.

For an infinite group, these are the only possibilities for the number of ends and these possibilities are mutually exclusive. It is known that a finitely generated group $G$ and a finite index subgroup $H$ of $G$ have the same number of ends, since the number of ends is a quasi-isometry invariant.  

We now state standard facts about amalgamated free products and HNN extensions that correspond to Theorem \ref{relatively-malnormal} and Lemma \ref{abelian-lemma}.  
We begin with the virtual malnormality results, corresponding to Theorem \ref{relatively-malnormal}. These results easily follow from the existence of normal forms
for amalgamated free products and HNN extensions, see Chapter IV of Lyndon and Schupp \cite{lyndon-schupp}.

\begin{lemma} Let $G$ be a finitely generated group that admits a non-trivial amalgamated free product decomposition $G =A\ast_C B$, where $C$ is finite. 
\begin{enumerate} 
\item If $g$, $h\in G$, then $gAg^{-1}\cap hBh^{-1}$ is finite;
\item If $g\in G\setminus A$, then $A \cap gAg^{-1}$ is trivial (and if $h\in G\setminus B$, then $B \cap hBh^{-1}$ is trivial).
\end{enumerate}
\label{afp-relatively-malnormal}
\end{lemma}

As a consequence, if $g\in G$ has infinite order and if $g^n\in A$ for some $n\ge 2$, then $g\in A$ (and similarly for $h\in B$).

\begin{lemma} Let $G$ be a finitely generated group that admits a non-trivial HNN decomposition $G =A\ast_C$, where $C$ is finite.
 If $g\in G\setminus A$, then $A\cap gAg^{-1}$ is trivial.
\label{hnn-relatively-malnormal}
\end{lemma}

As a consequence, if $g\in G$ has infinite order and if $g^n\in A$ for some $n\ge 2$, then $g\in A$.

We now state the lemmas corresponding to Lemma \ref{abelian-lemma}. The first follows from H. Neumann's generalization of the Kurosh theorem for amalgamated free products, see Lyndon and Schupp \cite{lyndon-schupp}, Chapter I, Proposition 11.22. The second follows from Britton's Lemma given in Lyndon and Schupp \cite{lyndon-schupp}, Chapter IV, Section 2.

\begin{lemma} Let $G$ be a finitely generated group that admits a non-trivial amalgamated free product decomposition $G =A\ast_C B$, where $C$ is finite. 
 If $K$ is an abelian subgroup of $G$ of rank at least two, then $K$ is contained in a conjugate of either $A$ or $B$.
\label{afp-abelian-lemma}
\end{lemma}

\begin{lemma} Let $G$ be a finitely generated group that admits a non-trivial HNN decomposition $G =A\ast_C$, where $C$ is finite.  If $K$ is an abelian subgroup
 of $G$ of rank at least two, then $K$ is contained in a conjugate of $A$.
\label{hnn-abelian-lemma}
\end{lemma}

We are now able to state and prove the analogue of Theorem \ref{our-main-theorem} for the one-endedness of groups with connected commutativity graph.  
This proof follows very much the same line of argument as the proof of Theorem \ref{our-main-theorem}.

\begin{theorem} Let $G$ be a finitely generated group which is not virtually ${\bf Z}$. Suppose there exists a (possibly infinite) generating set $S$ of 
cardinality at least two such that every element of $S$ has infinite order and $K(G, S)$ is connected. Suppose further that there exist adjacent vertices 
$s$, $s'$ of $K(G,S)$ and non-zero integers $n_s$, $n_{s'}$ so that $\langle s^{n_s}, (s')^{n_{s'}}\rangle$ is rank $2$ abelian.  
Then, $G$ has one end. 
\label{main-end-theorem}
\end{theorem}

\begin{proof} [of Theorem \ref{main-end-theorem}]
Since $G$ is assumed not to be virtually ${\bf Z}$, either $G$ has one end or $G$ has infinitely many ends.  Suppose for contradiction that 
$G$ has infinitely many ends, so that either $G$ admits a non-trivial amalgamated free product decomposition $G =A\ast_C B$ with $C$ finite, or 
$G$ admits a non-trivial HNN extension $G =A \ast_C$ with $C$ finite.  We give full details for the amalgamated free product case; the details in the 
HNN extension case are analogous.

We first show that no conjugate of either $A$ (or of $B$, the details are the same) can contain a non-zero power of an element of $S$. So, suppose 
there is some $g\in G$, some $s_0\in S$, and some $k\ne 0$ so that $s_{0}^{k} \in g A g^{-1}$.  (Note that, by the comment immediately following Lemma  
\ref{afp-relatively-malnormal}, this implies that $s_0\in gA g^{-1}$ as well.)  

Let $s_{1}$ be any vertex of $K(G, S)$ adjacent to
$s_{0}$. As there are non-zero integers $n_0$, $n_1$ so that $\langle s_{0}^{n_0}, s_{1}^{n_1}\rangle$ is abelian, we see that 
$\langle s_{0}^{n_0 k}\rangle \subseteq g A g^{-1} \cap s_{1}^{n_1}g A g^{-1} s_{1}^{-n_1} = g( A \cap g^{-1}s_{1}^{n_1}g A g^{-1} s_{1}^{-n_1}g) g^{-1}$. However, since 
$s_{0}$ has infinite order, the second part of Lemma \ref{afp-relatively-malnormal} implies $g^{-1} s_1^{n_1} g \in A$.  
By the comment immediately following Lemma \ref{afp-relatively-malnormal}, we see that $s_1\in gA g^{-1}$ as well.

Let $s$ be any element of $S$. By the connectivity of $K(G, S)$,
there is a sequence of elements $s_{0}, s_{1}, \ldots, s_{n}=s$
of $S$ such that $s_{k-1}$ and $s_k$ are adjacent in $K(G,S)$ for each $1\le k\le n$. The argument we have given above implies that if $s_{k-1} \in g A g^{-1}$,
then $s_{k} \in g A g^{-1}$. In particular, we have that $s\in g A g^{-1}$. Since $G$ is generated by $S$, it follows that $G$
and $g A g^{-1}$ are equal, contradicting the fact that
the subgroup $A$ is proper.  We conclude that if $G$ admits a non-trivial amalgamated free product decomposition $G =A\ast_C B$ with $C$ finite, then no conjugate of either $A$ or $B$ can contain a non-zero power
of an element of $S$. 

Now, by assumption there exist adjacent vertices $t$ and $t'$ of $K(G,S)$ for which there exist non-zero integers $n_t$, $n_{t'}$ so that $D =\langle t^{n_t}, (t')^{n_{t'}}\rangle$ is rank $2$ abelian.  By Lemma \ref{afp-abelian-lemma}, we see that $D$ must lie in some conjugate of  $A$ (or of $B$). In particular, this conjugate of $A$ contains a non-zero power of an element of $S$, and we have a final contradiction.
\end{proof}

\section{Applications}

In this section, we apply Theorems \ref{our-main-theorem} and \ref{main-end-theorem} to a selection of finitely generated groups, and deduce that  each is not strongly hyperbolic relative to any finite
collection of proper finitely generated subgroups (we will just say that such a group is {\em not strongly relatively hyperbolic}) and has one end. 

\subsection{Mapping class groups}
\label{mcg}

Let $\Sigma$ be a connected, orientable surface without
boundary, of finite topological type and negative Euler
characteristic.  As such, $\Sigma$ is the complement in a closed, orientable surface of a (possibly empty) finite set of points. The {\em mapping class group} ${\rm MCG}(\Sigma)$ associated to $\Sigma$
is the group of all homotopy classes of orientation preserving
self-homeomorphisms of $\Sigma$. For a thorough account of these groups, we refer the reader to Ivanov \cite{ivanov-map}. It is known that every mapping
class group ${\rm MCG}(\Sigma)$ is finitely presentable and can be
generated by Dehn twists. Masur and Minsky \cite{masur-minsky-i} prove that ${\rm MCG}(\Sigma)$ is weakly hyperbolic relative to a finite 
collection of curve stabilisers. 



Now let $S$ be the collection of primitive Dehn
twists about all elements of $\pi_1(\Sigma)$ that are represented by simple closed curves on $\Sigma$.  (Here, an element of a group is {\em primitive} if 
it is not a proper power of another element of the group.)  The associated commutativity graph is precisely the
$1$-skeleton of the {\em curve complex}, introduced in Harvey
\cite{harvey-cc}: this follows from the observation that two Dehn twists commute if and
only if their associated curves are disjoint. Moreover, the Dehn twists associated to any pair of adjacent vertices in the curve complex generate a rank
 $2$ free abelian group.  Such a graph is
connected provided $\Sigma$ is not a once-punctured torus or a
four-times punctured sphere. Hence, we have the following:

\begin{proposition} Let $\Sigma$ be a connected, orientable surface without boundary, of finite topological type and negative Euler characteristic.  If $\Sigma$ is not a once-punctured torus or a four-times punctured sphere, then the mapping class group ${\rm MCG}(\Sigma)$ of $\Sigma$ is not strongly relatively hyperbolic.
\label{mcg-rel-hyp}
\end{proposition}

This answers Question 6.24 of Behrstock \cite{behrstock} in the negative.
Note that, when $\Sigma$ is a once-punctured torus or a four-times punctured sphere, its mapping class group is isomorphic to ${\rm PSL}(2,{\bf{Z}})$ which
is a hyperbolic group. We remark that the result of Proposition \ref{mcg-rel-hyp} was previously obtained by Bowditch \cite{bowditch-hyperbolic}, using arguments
based on convergence groups. 

Since the mapping class group of a punctured sphere can be viewed as a braid group, the braid group $B_{n}$ on $n$ 
strings is not strongly relatively hyperbolic whenever $n \geq 5$.  This also follows by considering the usual presentation for $B_{n}$ and 
its corresponding commutativity graph.

We also have a corresponding result about ends.

\begin{proposition} Let $\Sigma$ be a connected, orientable surface without boundary, of finite topological type and negative Euler characteristic. 
 If $\Sigma$ is not a once-punctured torus or a four-times punctured sphere, then the mapping class group ${\rm MCG}(\Sigma)$ of $\Sigma$ has one end.
\label{mcg-end}
\end{proposition}

We note that Proposition \ref{mcg-end} is implicit in the work of Harer \cite{harer}, as it is proven there that ${\rm MCG}(\Sigma)$ is a virtual 
duality group with  virtual cohomological dimension greater than one, and such groups are known to have one end by a standard yoga.  It also is contained in  Culler and Vogtmann \cite{culler-vogtmann}.

As with Proposition \ref{mcg-rel-hyp}, the cases of the once-punctured torus and the four-times punctured sphere are anomalous; 
as noted above, in both of these cases, the mapping class group is isomorphic to ${\rm PSL}(2,{\bf{Z}})$, which is a free product and as such has infinitely many ends.



\subsection{The Torelli group}
\label{torelli}

The Torelli group $\mathcal{I}(\Sigma)$ of a connected, orientable surface $\Sigma$ is the kernel of the
natural action of the mapping class group ${\rm MCG}(\Sigma)$ on the first homology
group $H_{1}(\Sigma, {\bf Z})$. It is of continued interest, given
its connections with homology $3$-spheres and the number of basic
open questions it carries.  If $\Sigma$ is compact and has genus at least 3,  $\mathcal{I}(\Sigma)$ is generated by all Dehn twists around separating
 simple closed curves and all double twists around pairs of disjoint simple closed nonseparating curves
(called {\em bounding pairs}) that together separate (see \cite{johnson}). 

Farb and Ivanov \cite{farb-ivanov-torelli} introduce a graph they call the {\em Torelli geometry}. The vertices of this
graph comprise all separating curves and bounding pairs in
$\Sigma$, with two distinct vertices declared adjacent if their
corresponding curves or bounding pairs are disjoint. Whenever
$\Sigma$ has genus at least three this graph is
connected (this holds even when $\Sigma$ has non-empty boundary, see \cite{masur-schleimer}). For this reason, let us take $S$ to be the collection of primitive Dehn twists about separating curves and double twists
 around bounding pairs. The
corresponding commutativity graph $K(\mathcal{I}(\Sigma), S)$ is
precisely the Torelli geometry.  Also, as is the case with mapping class groups, adjacent vertices generate a rank 2 free abelian subgroup of
 $\mathcal{I}(\Sigma)$.  Thus, we have:

\begin{proposition} If $\Sigma$ is a closed, orientable surface of genus at least three, then the Torelli group $\mathcal{I}(\Sigma)$ of $\Sigma$ is
 not strongly relatively hyperbolic.
\end{proposition}

We conjecture this extends to all surfaces $\Sigma$ of genus $g$ and $n$ punctures with $2g + n - 4 \geq 1$ (to exclude small surfaces).   
For this, one would need to establish the connectivity of the Torelli geometry.

We have also the following result on ends of the Torelli group. Since $\mathcal{I}(\Sigma)$ has infinite index in ${\rm MCG}(\Sigma)$, the fact that $\mathcal{I}(\Sigma)$
has one end is independent of Section \ref{mcg}.

\begin{proposition} If $\Sigma$ is a closed, orientable surface of genus at least three, then the Torelli group $\mathcal{I}(\Sigma)$ of $\Sigma$ has one end.
\end{proposition}

\subsection{The special automorphism group of a free group}

In this subsection, we use the notation and basic results from Gersten \cite{gersten} (without further reference). 
Let ${\bf{F}}_n$ be the free group on $n$ generators and consider
the automorphism group ${\rm{Aut}}({\bf{F}}_n)$ of ${\bf{F}}_n$.
Abelianisation gives a surjective homomorphism 
\[ {\rm{Aut}}({\bf{F}}_n) \rightarrow {\rm{Aut}}({\bf{Z}}^n) = {\rm GL}(n,{\bf{Z}}). \]
Composing 
with the sign of determinant map 
\[ {\rm GL}(n,{\bf{Z}})\rightarrow \{\pm 1 \}, \]
we obtain a surjective homomorphism 
\[ \varphi: {\rm{Aut}}({\bf{F}}_n) \rightarrow \{\pm 1 \}, \]
which we call the {\em{determinant map}}.  

The {\em{special automorphism group}} of ${\bf{F}}_n$ is
${\rm{Aut}}^+({\bf{F}}_n) = {\rm{ker (\varphi)}}$, and has the following finite
presentation  in terms of
{\em{Nielsen maps}}: Let $X$ be a free basis
for ${\bf{F}}_n$ and let $E= X \cup X^{-1}$. Given $a$, $b \in E$
with $a \neq b$, $b^{-1}$, define the {\em{Nielsen map}} $E_{ab}$ for $a$, $b$
by $E_{ab}: {\bf{F}}_n \rightarrow {\bf{F}}_n$, where $E_{ab}(a)= ab$ and $E_{ab}(c) = c$ for $c \neq a$, $a^{-1}$.

Gersten \cite{gersten} shows that ${\rm{Aut}}^+({\bf{F}}_n)$
is generated by the finite set 
\[ S = \{E_{ab} \mid a, b \in E\mbox{ with }a \neq b, b^{-1}  \} \]
and that the following relation holds:
\[ [E_{ab}, E_{cd}] = 1 \mbox{ if }a \neq c, d, d^{-1}\mbox{ and }b \neq c, c^{-1}. \] 
(We suppress the full set of relations). 

Note that if $E_{ab}, E_{cd}$ are distinct and commute, they generate a rank 2 abelian subgroup of ${\rm{Aut}}^+({\bf{F}}_n)$, since
both of them have infinite order and neither is a power of the other. We then have:


\begin{proposition} With $S$ as above, $K({\rm{Aut}}^+({\bf{F}}_n), S)$ is connected
for $n \geq 5$. In particular, ${\rm{Aut}}^+({\bf{F}}_n)$ is not strongly relatively hyperbolic for $n\ge 5$.
\label{automorphism}
\end{proposition}

\begin{proof}
Let $E_{ab}$ and $E_{cd}$ be any two vertices of
$K({\rm{Aut}}^+({\bf{F}}_n), S)$. We have that $E_{ab}$ and
$E_{cd}$ are adjacent in $K({\rm{Aut}}^+({\bf{F}}_n), S)$ unless
$a \in \{ c, d, d^{-1} \}$ or $b \in \{c, c^{-1}\}$. Let us
consider $a=c$ (the remaining cases are similar). We
want to find a path from $E_{ab}$ to $E_{ad}$ in
$K({\rm{Aut}}^+({\bf{F}}_n), S)$. Since $n \geq 5$, there are $e,f
\in E \setminus \{a^{\pm 1}, b^{\pm 1}, d^{\pm 1} \}$. It then
follows, from the commutativity relation in
${\rm{Aut}}^+({\bf{F}}_n)$ mentioned above, that the sequence of
generators $E_{ab}, E_{ed}, E_{bf}, E_{ad}$ gives a path in $K({\rm{Aut}}^+({\bf{F}}_n), S)$ from
$E_{ab}$ to $E_{ad}$.
\end{proof}

Let ${\rm{Out}}^+({\bf{F}}_n) = {\rm{Aut}}^+({\bf{F}}_n) /
{\rm{Inn}}({\bf{F}}_n)$ be the {\em{special outer automorphism
group of ${\bf{F}}_n$}}. It is immediate that the natural
surjective homomorphism ${\rm{Aut}}^+({\bf{F}}_n)$ to
${\rm{Out}}^+({\bf{F}}_n)$ preserves the connectivity of our
commutativity graph for ${\rm Aut}^{+}({\bf{F}}_{n})$ and the
hypotheses of Theorem \ref{our-main-theorem}. We therefore deduce
the following:

\begin{corollary}
${\rm{Out}}^+({\bf{F}}_n)$ is not strongly relatively hyperbolic for $n \geq 5$.
\end{corollary}

Restricting the surjective homomorphism ${\rm{Aut}}({\bf{F}}_n)
\rightarrow {\rm GL}(n,{\bf{Z}})$ to ${\rm{Aut}}^+({\bf{F}}_n)$, we
obtain a homomorphism ${\rm{Aut}}^+({\bf{F}}_n) \rightarrow {\rm SL}(n,{\bf{Z}})$. The generating set
$S = \{E_{ab} \mid a, b \in E\mbox{ with }a \neq b, b^{-1}  \}$ projects onto a generating set $\overline{S}$ for $ {\rm SL}(n,{\bf{Z}})$ whose elements have infinite order in ${\rm SL}(n,{\bf{Z}})$, as immediately follows from the definition of 
the Nielsen maps. Also $K({\rm SL}(n,{\bf{Z}}), {\overline{S}})$ is connected, since $K({\rm{Aut}}^+({\bf{F}}_n), S)$ is.
Thus we have:

\begin{corollary} ${\rm SL}(n,{\bf{Z}})$ is not strongly relatively hyperbolic for $n \geq 5$.
\end{corollary}


Similarly, we have the following result about the number of ends of these groups.  We can extend the discussion to the 
automorphism and outer automorphism groups of ${\bf{F}}_n$ in this case, since the number of ends of a group is invariant under passing
to finite index subgroups.  We note that these results were earlier obtained by Culler and Vogtmann \cite{culler-vogtmann}.

\begin{proposition} For $n\ge 5$, the groups ${\rm{Aut}}({\bf{F}}_n)$, ${\rm{Aut}}^+({\bf{F}}_n)$, ${\rm{Out}}({\bf{F}}_n)$, ${\rm{Out}}^+({\bf{F}}_n)$, and ${\rm SL}(n,{\bf{Z}})$ are one-ended. 
\end{proposition}






\subsection{The Heisenberg group}

Recall that the {\em{3-dimensional Heisenberg group}} $\cal{H}$ is
given by the presentation 
\[ {\cal H} = \langle a,b,c \mid [a,b] = c,
[a,c]=1=[b,c] \rangle. \]
Consider the generating set $S= \{a,b,c\}$. It is evident from this presentation that the
commutativity graph $K({\cal{H}}, S)$ is connected.  Also, the group $\langle a, c\rangle$ is rank 2 abelian.

\begin{proposition}
The $3$-dimensional Heisenberg group $\cal{H}$ is not strongly
relatively hyperbolic.
\end{proposition}

We also have the following, which was previously known, see for instance Apurara \cite{apurara}. 

\begin{proposition}
The $3$-dimensional Heisenberg group $\cal{H}$ has one end.
\end{proposition}

\subsection{Thompson's group $F$}

Thompson's group $F$ is a torsion-free group of orientation-preserving, piecewise-linear homeomorphisms of the
unit interval of the real line; see Brown and Geoghegan \cite{thompson} or Cannon, Floyd and Parry \cite{cfp} for a complete
definition. Even though $F$ is a finitely presented group, it is sometimes convenient to work with the
infinite presentation 
\[ F = \langle x_0, x_1, x_2, \ldots \mid x_{j+1} = x_ix_jx_i^{-1}, i<j \rangle. \]
Let $S = \{x_j\mid j \ge 0 \}$. Clearly, the commutativity graph $K(F,S)$ is far from
being connected. So, we tinker: If we consider the generating set $S' =
S \cup \{x_0x_1^{-1}\}$, then the commutativity graph
$K(F,S')$ is still not connected, as $x_0$ and $x_1$ are isolated vertices.  However, $K(F,S')
\setminus \{x_0,x_1\}$ is connected, since $x_0x_1^{-1}$ commutes
with $x_i$ for all $i \geq 2$ (see, for instance, Burillo \cite{burillo}). This is enough for us to modify our main argument and deduce:


\begin{proposition}
Thompson's group $F$ is not strongly relatively hyperbolic.
\end{proposition}

\begin{proof}
Suppose, for contradiction, that $F$ were strongly hyperbolic
relative to the proper finitely generated subgroups $L_1, \ldots, L_p$. Since all
abelian subgroups of rank at least 2 are conjugate into some
$L_m$ by Lemma \ref{abelian-lemma}, and since $\langle x_0x_1^{-1}, x_2
\rangle$ is abelian of rank 2, we find that $\langle x_0x_1^{-1}, x_2 \rangle
\subset g L_m g^{-1}$ for some $m=1, \ldots, p$ and some $g\in F$. For
all $j\geq 2$, $x_0x_1^{-1}$ and $x_j$ commute and so
$\langle x_0x_1^{-1} \rangle \subset g L_m g^{-1} \cap x_j gL_m g^{-1} x_j^{-1}$.
Since $\langle x_0x_1^{-1}\rangle$ is infinite,  we have $x_j \in
gL_m g^{-1}$, for all $j\geq 2$, by Theorem \ref{relatively-malnormal}.

Now, $gL_m g^{-1}$ cannot contain both $x_0$ and $x_1$, since $S =\{x_j\mid
j\geq 0\}$ generates $F$ and $gL_m g^{-1}$ is a proper subgroup of $F$ by
assumption. So suppose $x_0 \notin g L_m g^{-1}$. (The case $x_1 \notin
g L_m g^{-1}$ is similar.) From the presentation above, we see
that $x_{j+1} = x_0 x_j x_0^{-1}$ for all $j\geq 2$. Therefore
$x_{j+1} \in x_0 g L_m g^{-1}x_0^{-1}$, since we have shown that $x_j \in
g L_m g^{-1}$, for all $j\geq 2$. Therefore $g L_m g^{-1} \cap x_0 g L_m g^{-1} x_0^{-1}$ is
infinite (as it contains $\langle x_j\rangle$ for any $j\ge 2$), contradicting Theorem \ref{relatively-malnormal}.
\end{proof}

Using the same style of argument, we also have the following.

\begin{proposition}
Thompson's group $F$ has one end.
\end{proposition}

\begin{proof}
Suppose, for contradiction, that $F$ has more than one end.  Since $F$ is not virtually ${\bf Z}$ and $F$ is torsion free, we see that $F$ admits a non-trivial free product splitting $F =A\ast B$.   Since all
abelian subgroups of rank at least 2 are conjugate into either $A$ or $B$ by Lemma \ref{afp-abelian-lemma}, and since $\langle x_0x_1^{-1}, x_2
\rangle$ is abelian of rank 2, we find that $\langle x_0x_1^{-1}, x_2 \rangle
\subset g A g^{-1}$ for some $g\in F$. (The details are similar if $\langle x_0x_1^{-1}, x_2 \rangle
\subset g B g^{-1}$  for some $g\in F$.)  For
all $j\geq 2$, $x_0x_1^{-1}$ and $x_j$ commute and so
$\langle x_0x_1^{-1} \rangle \subset g A g^{-1} \cap x_j gA g^{-1} x_j^{-1}$.
Since $\langle x_0x_1^{-1}\rangle$ is infinite,  we have $x_j \in
gA g^{-1}$, for all $j\geq 2$, by Lemma \ref{afp-relatively-malnormal}.

Now, $gA g^{-1}$ cannot contain both $x_0$ and $x_1$, since $S =\{x_j\mid
j\geq 0\}$ generates $F$ and $gA g^{-1}$ is a proper subgroup of $F$ by
assumption. So suppose $x_0 \notin g A g^{-1}$. (The case $x_1 \notin
g A g^{-1}$ is similar.) From the presentation above, we see
that $x_{j+1} = x_0 x_j x_0^{-1}$ for all $j\geq 2$. Therefore
$x_{j+1} \in x_0 g A g^{-1}x_0^{-1}$, since we have shown that $x_j \in
g A g^{-1}$, for all $j\geq 2$. Therefore $g A g^{-1} \cap x_0 g A g^{-1} x_0^{-1}$ is
infinite (as it contains $\langle x_j\rangle$ for any $j\ge 3$), contradicting Lemma \ref{afp-relatively-malnormal}.
\end{proof}

{\footnotesize J. W. Anderson (corresponding author)\\School of Mathematics\\University of Southampton\\Southampton SO17 1BJ\\England\\j.w.anderson@soton.ac.uk}

{\footnotesize J. Aramayona\\Mathematics Institute\\University of Warwick\\Coventry CV4 7AL\\England\\jaram@maths.warwick.ac.uk}

{\footnotesize K. J. Shackleton\\School of Mathematics\\University of Southampton\\Southampton SO17 1BJ\\England\\k.j.shackleton@maths.soton.ac.uk}

\end{document}